\documentclass[12pt]{article}
\usepackage{amsmath, amsthm, amssymb,  graphicx}

\usepackage{geometry} 

\title{Continuous Interval Exchange Actions}
\author{Christopher F. Novak}

\newtheorem{theorem}{Theorem}[section]

\newtheorem{lemma}[theorem]{Lemma}
\newtheorem{proposition}[theorem]{Proposition}
\newtheorem{Conjecture}[theorem]{Conjecture}

\theoremstyle{remark}

\theoremstyle{definition}

\def\co{\colon\thinspace}

\begin{document}

\maketitle

\begin{abstract}    
\noindent Let $\mathcal{E}$ denote the group of all interval exchange transformations on $[0,1)$. Given a suitable topological group structure on 
$\mathcal{E}$, it is possible to classify all one-parameter interval exchange actions (continuous homomorphisms $\mathbb{R} \rightarrow \mathcal{E}$). In 
particular, up to conjugacy in $\mathcal{E}$, any one-parameter interval exchange action factors through a rotational torus action. 
\end{abstract}


\section{Introduction}

\begin{figure}[htb]   
\begin{center}
\includegraphics{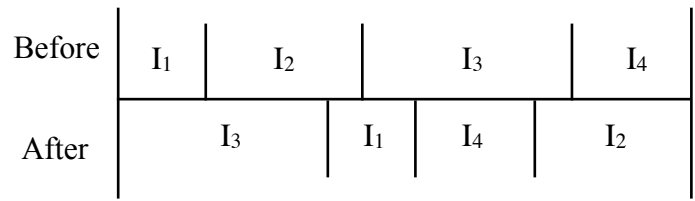}
\caption{An interval exchange with $\pi = (1\ 2\ 4\ 3)$   \label{fig interval exchange}  }
\end{center}
\end{figure}

An interval exchange transformation is a map $[0,1) \rightarrow [0,1)$   
defined by a finite partition of the unit interval into half-open intervals and 
 a rearrangement of these intervals by translation.        
See Figure~\ref{fig interval exchange} for a graphical example.

The dynamics of interval exchanges were first studied in the late seventies by
 Keane {\bf \cite{Keane75}} {\bf\cite{Keane77}}, Katok {\bf \cite{Katok80}}, Rauzy {\bf \cite{Rauzy79}}, Veech {\bf \cite{Veech78}}, and others. 
 This initial stage of research culminated in the independent proofs by Masur {\bf \cite{Masur82}} and Veech {\bf\cite{Veech82}} that
 almost every interval exchange is uniquely ergodic.
 See the recent survey of Viana {\bf \cite{Viana06}} for a unified presentation of these results. 
  
  There is currently much interest and activity in the dynamics of interval exchanges. This is due in part to the recent resolution of certain long-standing problems in this area; one important example is the work of Avila and Forni {\bf\cite{AvilaForni07}} in which they prove that almost every interval exchange is weakly mixing. Much of the study of interval exchanges is closely related to dynamics 
 on the moduli space of translation surfaces; an introduction to this topic and its connection to interval exchanges is found in a survey of
 Zorich {\bf\cite{Zorich06}}.  
 
An extension of the study of single interval exchanges is to consider their dynamics in terms of group actions.
The set $\mathcal{E}$ of all interval exchange transformations forms a group under composition, and an \emph{interval exchange 
action} of a group $G$ is a homomorphism $G \rightarrow \mathcal{E}$. A general and fundamental question is to determine if a given
group $G$ has faithful interval exchange actions. More broadly, it can be asked if there are
 general algebraic obstructions to the existence of such actions.  On the other hand, if such actions do exist for a given group, it is
  desirable to attempt to classify them in some way. The goal of this paper is to classify continuous 
  interval exchange actions of the group $\mathbb{R}.$

 The study of interval exchange actions is motivated by the study of other transformation groups, particularly groups of 
 homeomorphisms and diffeomorphisms of one-dimensional manifolds. However, what is known about the structure of $\mathcal{E}$ suggests
  that there may be substantial differences between $\mathcal{E}$ and these groups. For instance, it is shown
   in {\bf \cite{Novak09}} that no subgroup of $\mathcal{E}$ has distortion elements. In contrast, 
the groups $\text{Diff}^{\,\omega}(\mathbb{R})$ and $\text{Diff}^{\,\omega}(S^1)$ of real-analytic diffeomorphisms on the line and circle both
 contain such elements. See {\bf \cite{Franks06}} for definitions, examples, and results involving actions of groups having distortion elements.
  
 In addition, many basic questions that are well understood for diffeomorphisms of 1--manifolds
  are currently open for the group $\mathcal{E}$. For instance:
  
  \begin{itemize}
  \item[(1)] Does $\mathcal{E}$ contain a free subgroup on two generators? (Katok)
  \item[(2)] Does $\mathcal{E}$ contain groups of intermediate growth? (Grigorchuk)
  \item[(3)] Is every solvable subgroup of $\mathcal{E}$ virtually abelian? (Navas)
   \end{itemize}
   
  For question (1), it is easy to construct examples of non-abelian free groups in $\text{Diff}(S^1)$ or $\text{Diff}(\mathbb{R})$ by means of the ping-pong construction. 
 More detailed results, analogous to the Tits' alternative, are also known for $\text{Homeo}_+(S^1)$ and $\text{Diff}^{\, \omega}(S^1);$ see
  {\bf \cite{Margulis00}}  and {\bf \cite{FarbShalen02}}, respectively. Question (2) is answered in the affirmative for the group
  $\text{Diff}_+^{\,1}([0,1])$ in {\bf \cite{Navas08}}. This work also shows that for any $\alpha > 0,$ any subgroup of $\text{Diff}_+^{\,1+\alpha}([0,1])$ with sub-exponential growth must  be virtually nilpotent. This gives a negative answer to question (2) in this case, due to the fundamental result of Gromov {\bf \cite{Gromov81}} that the finitely generated virtually nilpotent groups are exactly those having polynomial growth. Question (3) is also well understood for various 
  transformation groups on 1--manifolds; for instance, see {\bf \cite{Bleak08}}, {\bf \cite{BurslemWilkinson04}}, and {\bf \cite{Navas04}}.

To introduce the results of the current work, consider the following precise definition and notation for an interval exchange.
Let $\pi \in \Sigma_n$ be a permutation of $\{1,2, \ldots, n\}$, and let $\lambda$ be a vector in the simplex 
\[ \Lambda_n = \left\lbrace \lambda = (\lambda_1, \ldots, \lambda_n)\co \lambda_i > 0, \
  \sum \lambda_i = 1     \right\rbrace \subseteq \mathbb{R}^n.  \]
   The vector $\lambda$ induces a partition of $[0,1)$ into intervals
  \begin{equation} \label{IET_form1} I_j = \left[ \beta_{j-1} := \sum_{i=1}^{i = j-1} \lambda_i,\ \, \beta_j := 
  \sum_{ i=1}^{i = j} \lambda_i \right), \ \ 1 \leq j \leq n. 
  \end{equation}  
Let $f_{(\pi, \lambda)}$ be the interval exchange that translates each $I_j$ such that the ordering of these intervals
within $[0,1)$ is permuted according to $\pi$. More precisely, 
\begin{equation} \label{IET_form3} f_{(\pi, \lambda)} (x) = x + \omega_j, \quad \text{if }x \in I_j, \end{equation} 
where 
\begin{equation} \label{IET_form2} \omega_j = \Omega_\pi(\lambda)_j\, = \sum_{i:\, \pi(i) < 
\pi(j)} \lambda_i\ -\ \sum_{i:\, i<j} \lambda_i.  \end{equation}
Note that $\Omega_\pi: \Lambda_n \rightarrow \mathbb{R}^n$ is a linear map depending only on $\pi$.

\begin{figure}[htb] \begin{center} 
\includegraphics{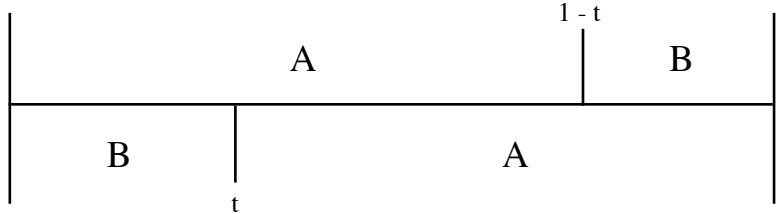}
\caption{The rotation $r_t$ with $\pi = (2\ 1)$ and $\lambda = (1-t, t)$    \label{fig basicrotation} }
\end{center}
\end{figure}

An initial example of an interval exchange $\mathbb{R}$-action is defined by 
\[     t \mapsto r_t,  \] 
where $r_t \in \mathcal{E}$ is the rotation $r_t(x) = x + t \text{ (mod 1)},$ as depicted in Figure~\ref{fig basicrotation}. This action is not 
faithful, but its kernel is a discrete subgroup of $\mathbb{R}$. 

To construct a faithful interval exchange action of $\mathbb{R}$, choose numbers $\alpha_1$ and $\alpha_2$ in $(0, 1)$ such that $\alpha_1/\alpha_2$ is irrational.
For a real number $t$, define $f_t$ as the map that rotates by $t\alpha_ 1$ (mod 1/2) on the interval $[0, 1/2)$ and rotates by $t\alpha_2$ (mod 1/2) 
on the interval $[1/2, 1)$. That is, $f_t = f_{(\pi,\, \lambda(t))}$ for the data
\[\pi = (1\ 2)(3\ 4), \quad 
		 \lambda(t) = \left( \frac{1 - \{ t\alpha_1\} }{2},\ \ \frac{\{ t\alpha_1\}}{2},\ \ \frac{1 - \{t\alpha_2 \}}{2},\ \  \frac{\{t\alpha_2\}}{2}   \right),  \]
 where $\{ \cdot \}$ denotes the fractional part of a real number.
It is easy to check that $t \mapsto f_t$ is a group homomorphism $\mathbb{R} \rightarrow \mathcal{E}$, and the action is faithful due to the assumption that 
$\alpha_1/\alpha_2$ is irrational.

\begin{figure}[htb]  
\begin{center}
\includegraphics{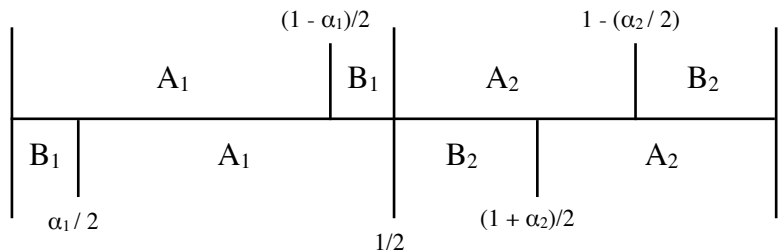}
\caption{The action $  (\alpha_1, \alpha_2) \mapsto  f_{(\pi,\, \lambda(\alpha_1, \alpha_2))}$   \label{torus action}  }
\end{center}
\end{figure}

This faithful $\mathbb{R}$-action can be viewed as the restriction of an action of the torus $\mathbb{T}^2 = \mathbb{R}^2 /\mathbb{Z}^2$ 
(identified with $[0,1) \times [0,1)$) defined by the map 
$ (\alpha_1, \alpha_2) \mapsto  f_{(\pi,\, \lambda(\alpha_1, \alpha_2))}, $
where 
\[ \pi = (1\ 2)(3\ 4), \quad    
  		\lambda(\alpha_1, \alpha_2) = \left(\frac{1 - \alpha_1}{2}, \ \frac{\alpha_1}{2},\  \frac{1- \alpha_2}{2},\  \frac{\alpha_2}{2}    \right).   \]
See Figure~\ref{torus action} for an illustration of the map 	$ f_{(\pi,\, \lambda(\alpha_1, \alpha_2))}$.	

A general class of torus actions can be defined by a similar construction. For any partition vector $\lambda \in \Lambda_n$ and for 
any $\alpha = (\alpha_1, \alpha_2, \ldots, \alpha_n) \in \mathbb{T}^n,$ define the interval exchange $f_{(\alpha, \lambda)}$ by 

\begin{equation} \label{torus action element} f_{(\alpha, \lambda)}\co x \mapsto \left\lbrace \begin{array}{ll}
 x+\lambda_j {\alpha}_j, & x\in [\beta_{j-1},\, \beta_j - \lambda_j{\alpha}_j ) \\          
 x+ \lambda_j{\alpha}_j - \lambda_j, 
 & x\in [\beta_j -\lambda_j{\alpha}_j,\, \beta_j ), \end{array} \right.    \end{equation}
where the points $\beta_j$ are the boundary points of the partition intervals defined by $\lambda$. The map $f_{(\alpha, \lambda)}$
is also defined by the permutation $\pi = (1\ 2)(3\ 4)\cdots((2n-1)\ 2n)$ and the partition vector 
\[ \left( \lambda_1(1 - \alpha_1),\, \alpha_1 \lambda_1,\,  \lambda_2(1 - \alpha_2),\, \alpha_2 \lambda_2,\, \ldots , 
 \lambda_n(1 - \alpha_n),\, \alpha_n \lambda_n           \right).  \] 

The action $\alpha \in \mathbb{T}^n \mapsto f_{(\alpha, \lambda)}$ is called the \emph {standard torus action} associated to $\lambda$. 
Restricting a standard torus action to a one-parameter subgroup gives an action of $\mathbb{R}$. 
A (one-parameter) \emph{rotation action} is defined as any action $\mathbb{R} \rightarrow \mathcal{E}$ that is conjugate in $\mathcal{E}$ to a one-parameter 
subgroup of a standard torus action. The image in $\mathcal{E}$ of a rotation action will be referred to as 
a rotation subgroup. 

The main result of this paper is that under some natural and unrestrictive topological assumptions,
 the rotation actions classify all continuous interval exchange actions of $\mathbb{R}$. To specify these conditions, define a permutation 
 $\pi \in \Sigma_n$ to be \emph{unpartitioned} if  $\pi(j+1) \neq \pi(j) + 1$, for all $j$ such that $1\leq j \leq n-1$. It is shown in Proposition~\ref{unique_coords} that an 
 interval exchange is defined by a unique pair $(\pi, \lambda)$ if one restricts to unpartitioned permutations. For each unpartitioned 
 $\pi \in \Sigma_n,$ there is a coordinate map  $\Gamma_\pi\co \Lambda_n \rightarrow \mathcal{E} $
defined by $\Gamma_\pi (\lambda) = f_{(\pi, \lambda)}$. The definition of the map $ f_{(\pi, \lambda)}$ in equations 
\eqref{IET_form1}--\eqref{IET_form2} 
extends to vectors $\lambda \in \overline{\Lambda}_n$ by allowing some of the partition intervals $I_i$ to be degenerate.  Thus, the coordinate 
maps extend to maps $\Gamma_\pi\co \overline{\Lambda}_n \rightarrow \mathcal{E}$. The needed topological conditions essentially require the 
$\Gamma_\pi$ to be continuous parameterizations.

\begin{theorem} \label{one-parameter IET subgroups}
Suppose that $\mathcal{E}$ has a topological group structure such that for every unpartitioned permutation $\pi \in \Sigma_n$, the coordinate map 
$\Gamma_\pi\co \overline{\Lambda}_n \rightarrow \mathcal{E}$ is continuous and the restriction  $\Gamma_\pi|_{\Lambda_n}$ is a homeomorphism onto its image. 
Then an action $\mathbb{R} \rightarrow \mathcal{E}$ is continuous if and only if it is a rotation action.
\end{theorem}

Due to this result, the image in $\mathcal{E}$ of a rotation action will be referred to as a one-parameter subgroup of $\mathcal{E}$.
Based on the classification in {\bf \cite{Novak09}} of interval exchange centralizers, it is not hard to see that if two one-parameter subgroups commute, 
then they are simultaneously conjugate to subgroups of a common standard torus action. Thus, the group generated by two distinct commuting one-parameter subgroups 
is conjugate to the image of a two-dimensional subgroup of $\mathbb{T}^n$ under some standard torus action. However, the situation for noncommuting
 one-parameter subgroups appears to be quite different. 
 		
\begin{Conjecture} \label{noncommuting rotation actions}
Let $F = \{f_t\}$ and $G = \{g_s\}$ be one-parameter subgroups of $\mathcal{E}$. If $F$ and $G$ do not commute, then the group $\langle F, G \rangle$ 
has elements that are not contained in any one-parameter subgroup of $\mathcal{E}$. 
\end{Conjecture}		

A consequence of this conjecture is the observation that up to conjugacy in $\mathcal{E}$, an interval exchange action of a 
connected Lie group must factor through a standard torus action.

A motivating example for this conjecture is the pair of one-parameter subgroups $F  = \{r_t\}$ and $G = \{r_{s, \delta}\},$ where the map $r_{s, \delta}$ denotes 
a restricted rotation by $s\delta\, (\text{mod }\delta)$ supported on the interval $[0, \delta)$. It can be shown that certain products in  $\langle F, G \rangle$, such as 
$h = r_t\circ r_{s, \delta}$ with sufficiently small $t$ and $s$, will have linear discontinuity growth; i.e., if $D(h^n)$ denotes the number of discontinuities of $h^n$, 
then $D(h^n) \sim Cn$ for some $C>0$. It is not difficult to see that interval exchanges with linear discontinuity growth  
cannot be in the image of a rotation action. Hence, by Theorem~\ref{one-parameter IET subgroups}, the group $\langle F, G \rangle$ contains elements that do not lie on any one-parameter interval exchange group.  

\subsubsection*{Acknowledgements}
The author is indebted to Professor John Franks for suggesting this topic of study and for many helpful conversations while completing this work. The author 
would like to thank the referee for many useful comments and for suggesting much important background information. 


\section{Coordinates and Topology of $\mathcal{E}$}

The restriction to unpartitioned permutations is suggested by the fact that such a permutation  $\pi$ properly describes the discontinuities of 
the map $f_{(\pi, \lambda)}$.

\begin{lemma} \label{unpart=discts} $\pi \in \Sigma_n$ is unpartitioned if and only if, for any $\lambda \in \Lambda_n$,
the interval exchange $f_{(\pi,\lambda)}$ is discontinuous, as a map $[0,1) \rightarrow [0,1)$, at
 precisely each of $\beta_1,\, \ldots, \beta_{n-1}$.   \end{lemma}

\noindent \emph{Proof}: If $\pi(j+1) = \pi(j)+1$ for some $j$, then for any $\lambda \in \Lambda_n,$
 the map $f_{(\pi, \lambda)}$ restricts to a translation on  $I_j \cup I_{j+1}.$ In particular, $f_{(\pi,\lambda)}$ is continuous at $\beta_j$.

Conversely,  if $f = f_{(\pi, \lambda)}$ is 
continuous at $\beta_j$, then both $I_j$ and $I_{j+1}$ are translated the same distance by $f$. Consequently, 
 $\pi(j+1) = \pi(j) + 1$, which implies $\pi$ is not unpartitioned.  $\square$

\begin{proposition} \label{unique_coords} For any interval exchange $f \in \mathcal{E}$, there exists a positive integer $n,$
 an unpartitioned $\pi \in \Sigma_n$, and $\lambda \in \Lambda_n$, all of which are unique, such
that $f = f_{(\pi, \lambda)}$.  \end{proposition}

\noindent \emph {Proof}: To show the existence of $n, \pi,$ and $\lambda$, let 
\[ 0 < \beta_1 < \beta_2 < \, \ldots\, < \beta_{n-1} < 1 \]
be the finite set of points in $(0,1)$ at which $f$ is discontinuous as a map $[0,1) \rightarrow [0,1);$
 this defines $n$. Setting $\beta_0 = 0$ and 
$\beta_n = 1$, define $\lambda \in \Lambda_n$ by 
\[ \lambda_j = \beta_j - \beta_{j-1}, \ \ \ j = 1, \, \ldots,n. \]

The permutation $\pi$ is defined by the reordering of the points $\beta_{i-1}$ induced by the map $f$. Thus, 
$\pi(i) = j$ if and only if $\#\{k\co f(\beta_k) < f(\beta_{i-1}) \} +1 = j$. By construction $f = f_{(\pi, \lambda)}, $ and
 $\pi$ is unpartitioned by Lemma~\ref{unpart=discts}, since $f$ is discontinuous at precisely $\beta_1,\, \ldots, \beta_{n-1}$.
The uniqueness of $\pi$ and $\lambda$ now follows since these data were constructed using intrinsic features of 
the transformation $f$.  \ \ $\square$ \\

The most natural choice of a topology on $\mathcal{E}$ that satisfies the conditions of Theorem~\ref{one-parameter IET subgroups}
 is the CW--complex topology induced by the cell maps $\Gamma_\pi$. 
It is not difficult to check that the group operations of $\mathcal{E}$ are continuous with respect to this topology. In particular, since the CW--complex structure
for $\mathcal{E}$ involves only countably many cells, the product $\mathcal{E} \times \mathcal{E}$ has a CW--complex structure induced by the cell maps $\Gamma_\pi \times \Gamma_\sigma$ (see {\bf \cite{Hatcher01}}, Prop A-2, p. 521). In addition, a proof of Theorem~\ref{one-parameter IET subgroups} with respect to 
the CW--complex topology is not difficult (particularly given Lemma~\ref{local_rotation} and Proposition~\ref{local_cts_orbits_iff_rotn} below); this is due to the fact 
that in the cell topology any compact subset must be contained within finitely many cells. 

\begin{figure}[htb]  
\begin{center}
\includegraphics{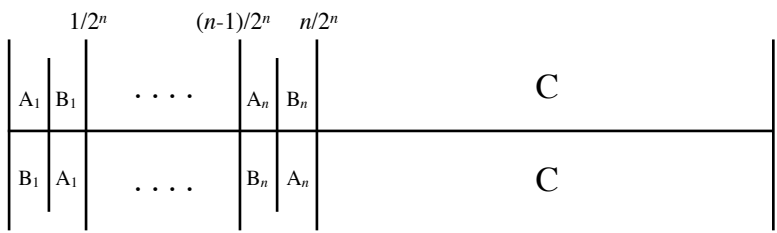}
\caption{The map $f_n$  \label{CW_counterexample}   }
\end{center}
\end{figure}

However, the CW topology on $\mathcal{E}$ does not properly reflect the action of $\mathcal{E}$ on various function spaces defined over $\mathbb{T}$. For instance, consider the following sequence $f_n$ of interval exchanges. 
 Define $\pi^{(1)} = (1\ 2)$ and $\lambda^{(1)} = (1/2, 1/2)$; let $f_1 = f_{(\pi^{(1)}, \lambda^{(1)})}$. For integers $n \geq 2$, let 
 \[ \pi^{(n)} = (1\ 2)\cdot(3\ 4)\cdots((2n-1)\ (2n)) \cdot(2n+1) \in \Sigma_{2n+1},  \]
\[ \lambda^{(n)} = \left( \frac{1}{2^{n+1}},\,  \frac{1}{2^{n+1}},\, \ldots ,\,  \frac{1}{2^{n+1}},\, 1 - \frac{n}{2^{n}} \right),  \]
and let $f_n = f_{(\pi^{(n)},\, \lambda^{(n)}) }$; see Figure~\ref{CW_counterexample}. 

Since each element of $\{f_n\}$ is in a different cell, this set is closed in the CW toplogy of $\mathcal{E}$. In particular, it does not contain the identity as a limit point, 
even though the mappings $f_n$ converge uniformly to the identity. 

\begin{figure}[htb] 
\begin{center}
\includegraphics{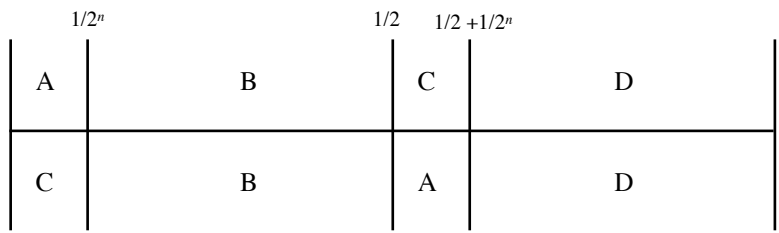}
\caption{The map $g_n$   \label{uniform_counterexample}  }
\end{center}
\end{figure}

The topology of uniform convergence is also not a suitable topology on $\mathcal{E}$ since the maps are not usually continuous. For instance, the sequence 
of maps $g_n = g_{(\sigma^{(n)}, \eta^{(n)})}$ defined by 
\[   \sigma^{(n)} = (1\ 3) \in \Sigma_4 , \quad \quad \eta^{(n)} = \left( \frac{1}{2^n}, \frac{2^{n-1}-1}{2^n}, \frac{1}{2^n}, \frac{2^{n-1}-1}{2^n}  \right)   \]
does not converge uniformly to the identity; see Figure~\ref{uniform_counterexample}.

We now define a topological group structure on $\mathcal{E}$ in which both of the above sequences converge to the identity. 
Let $\rho$ denote the shortest-path metric on the circle $\mathbb{T},$ identified with $[0,1)$. Given $f, g \in \mathcal{E}$, define 
\[  d(f, g) = \int_{\mathbb{T}} \rho (f(x), g(x)) \, d\mu(x),  \]
where $\mu$ denotes Lebesgue measure.

\begin{proposition} \label{E is topological group} The function $d$ is a metric on $\mathcal{E},$ and the metric space $(\mathcal{E}, d)$
 is a topological group.   \end {proposition}

\noindent \emph{Proof}: The nonnegativity, symmetry, and triangle inequality for $d$ follow from the corresponding properties of the metric $\rho$.
Moreover, note that $d(f,g) = 0$ implies that $f(x) = g(x)$ $\mu$--a.e., and two interval exchanges that coincide $\mu$--a.e.\hspace{-0.04in} must be identical. Thus $d$ 
is a metric.

To show that composition is continuous with respect to $d$, suppose there are convergent
 sequences $f_n \rightarrow f$ and $g_n \rightarrow g$. It suffices to estimate for all sufficiently large $n$ that $f_n(g_n(x))$ is close to $f(g(x))$ outside  
 of a set with small measure, since the $\mathbb{T}$--metric $\rho$ is bounded on the exceptional set.
 To achieve this estimate, two comparisons can be made. First, for sufficiently large $n,$ the $\mathbb{T}$--distance between $f_n(g_n(x))$ and $f(g_n(x))$ is 
 small outside of a set with small measure since $d(f_n, f)$ is close to zero. Next, the distance between $f(g_n(x))$ and $f(g(x))$ is small outside of a set with 
 small measure since $d(g_n, g)$ is close to zero and since the map $f$ is a translation on most sufficiently small intervals.  Thus, composition is continuous 
 with respect to the metric $d$.

It remains to show that inversion in $\mathcal{E}$ is continuous. First, note that the metric $\rho$ is invariant under
\emph{right} translation in the group, since all interval exchange transformations preserve Lebesgue measure. Thus, 
\[ d(f, id) = d(id, f^{-1}).\]
Consequently, if $f_n \rightarrow {id}$, then $f_n^{-1} \rightarrow {id}$. Thus inversion is continuous at the identity.
In general, if $f_n \rightarrow f$, the continuity of composition implies that $f^{-1}f_n \rightarrow id$. But then 
$f_n^{-1}f \rightarrow id$, and applying the continuity of composition again yields $f_n^{-1} \rightarrow f^{-1}$, as desired.
$\square$\\

 \begin{proposition} \label{Coord_maps} For any unpartitioned permutation $\pi$, the
  map $\Gamma_\pi\co \overline{\Lambda}_n \rightarrow \mathcal{E}$ is continuous with respect to the metric $d$. Consequently, the 
 restriction $\Gamma_\pi|_{\Lambda_n}$ is a homeomorphism onto its image.      \end{proposition} 
 
 \noindent  \emph{Proof}:  
 It has been shown in Proposition~\ref{unique_coords} that the restriction $\Gamma_\pi|_{\Lambda_n}$ is injective. If the map
 $\Gamma_\pi$ is continuous, then the restriction  $\Gamma_\pi|_{\Lambda_n}$ is a homeomorphism onto its image due to the compactness of 
   $\overline{\Lambda}_n$.
 
 To show the continuity of $\Gamma_\pi\co \overline{\Lambda}_n \rightarrow \mathcal{E}$, suppose
 that $\lambda^{(m)} \rightarrow \lambda$ in $\overline{\Lambda}_n$, and let $f^{(m)}$ and $f$ denote 
 $\Gamma_\pi(\lambda^{(m)})$ and $\Gamma_\pi(\lambda)$, respectively. 
 Given some $\epsilon >0$, for all sufficiently large $m,$ we have
\[ \left| \lambda_j - \lambda^{(m)}_j \right| < \frac{\epsilon}{n},\ \ j=1, \ldots, n. \]
Comparing the difference between boundary points of the partition intervals of $f^{(m)}$ and
$f$, we have
\[  \left| \beta_j - \beta_j^{(m)} \right| = \left| \sum_{k=1}^j \lambda_k - \sum_{k=1}^j \lambda^{(m)}_k   \right|
  \leq  \sum_{k=1}^j \left| \lambda_k - \lambda^{(m)}_k  \right| < \epsilon.    \]
Thus, for sufficiently large $m,$ the partition intervals
$I_j$ and $I^{(m)}_j$ overlap up to a set of small measure. That is,
\[ \mu\left(I_j \setminus I^{(m)}_j\right) < 2\epsilon. \]
 
Next, observe that the translation vectors $\omega^{(m)} = \Omega_\pi(\lambda^{(m)})$ converge to 
$\omega = \Omega_\pi(\lambda)$, since the map $\Omega_\pi$ is linear. Thus, for all sufficiently large $m$, 
\[ \left| \omega_j - \omega^{(m)}_j   \right| < \epsilon. \]
Therefore,
\[ d(f,f^{(m)}) =  \sum_{j=1}^n \left( \int_{I_j \cap I^{(m)}_j}\rho(fx, f^{(m)}x)\, d\mu(x) 
 +  \int_{I_j \setminus I^{(m)}_j}\rho(fx, f^{(m)}x)\, d\mu(x)  \right) = \]

\[ \sum_{j=1}^n \int_{I_j \cap I^{(m)}_j}\rho(x + w_j, x + w^{(m)}_j)\, d\mu(x)
 + \sum_{j=1}^n \int_{I_j \setminus I^{(m)}_j}\rho(fx, f^{(m)}x)\, d\mu(x). \]
The first term in this last expression is bounded by $\epsilon$, since $\rho( x + w_j, x + w^{(m)}_j) < \epsilon$
on the sets $I_j \cap I^{(m)}_j$. The second term is bounded by $n\epsilon$, since $\rho \leq 1/2$ and
$\mu(I_j \setminus I^{(m)}_j) < 2\epsilon.$ Thus $d(f, f^{(m)}) < (n+1)\epsilon$ for all sufficiently large $m$, which
proves that $\Gamma_\pi$ is continuous. \ \  $\square$\\

In addition to having the desired topological properties with respect to the coordinate maps $\Gamma_\pi$, the topology induced
by the metric $d$ has another natural interpretation. Since all interval exchanges preserve Lebesgue measure on $[0,1)$, 
one may view the group $\mathcal{E}$ as a group of unitary operators on $L^2(\mathbb{T}, \mu)$.

\begin{proposition}  The topology induced on $\mathcal{E}$ by the metric $d$ coincides with the strong operator topology
when $\mathcal{E}$ is viewed as a subgroup of $B(L^2(\mathbb{T}, \mu))$.   \end{proposition}

\noindent \emph{Proof}: Let $\{f_n\}$ be a sequence in $\mathcal{E}$ and let $\{T_n\}$ be the corresponding sequence of unitary 
operators; similarly, let $f \in \mathcal{E}$ and let $T$ denote its corresponding operator. 

First, suppose that $d(f_n, f)$ converges to zero, and let $\phi \in L^2(\mathbb{T})$ be any continuous function. As in the above arguments, 
having $d(f_n, f)$ close to zero means that $\rho(f_n(x), f(x)) <\delta$ on a set $A_n$ whose complement has measure approaching zero as $n$ increases.
 Using the uniform continuity 
 of $\phi,$ it follows that $|T_n\phi(x) - T\phi(x)| < \epsilon$ on $A_n$. 
 Since $\phi$ is bounded by a constant $M$, the contribution to $\Vert T_n\phi - T\phi \Vert_2^2$
 from points in the complement of $A_n$ is bounded by $4M^2 \mu(A_n^c)$, which approaches zero. Thus $T^n$ converges to $T$ in the strong operator topology.  

Conversely, suppose  $T^n$ converges strongly to $T$. Given $\epsilon >0$, partition $\mathbb{T}$ into intervals $I_1, I_2, \ldots, I_M$ 
such that $f$ is continuous on each interval $I_j$ and $\mu(I_j) < \epsilon$ for all $j$. Since $\Vert T_n \chi_{I_j} - T \chi_{I_j} \Vert_2 \rightarrow 0,$
there exists $N_j$ such that
 \[ \mu(\{ x\in I_j: f_n(x) \in f(I_j) \}) > (1-\epsilon)\mu(I_j)\]
  for all $n\geq N_j$. Since $f(I_j)$ is an interval, we have $\rho(f_n(x), f(x)) < \epsilon$ for any $x \in I_j$ such that  $f_n(x) \in f(I_j)$. 
    Let $N = \max\{N_j\}$. Then for any $n > N$, the set of points $x \in \mathbb{T}$ for which $\rho(f_n(x), f(x)) < \epsilon$ has measure greater than $(1-\epsilon).$
    Thus $d(f_n, f) < \epsilon + (1/2)\epsilon$, which implies $d(f_n, f) \rightarrow 0$. \quad $\square$


\section{Proof of the Classification Theorem}

In preparation for proving Theorem~\ref{one-parameter IET subgroups},
 it is useful to describe the dynamics of a rotation subgroup $F = \{f_t\}$
on $[0,1)$. Suppose $f_t = f_{([t\alpha], \lambda)}$, where $\alpha \in \mathbb{R}^n$ and $[t\alpha]$ is the equivalence class
 of $t\alpha$ in $\mathbb{T}^n \cong [0,1)^n.$ See equation \eqref{torus action element} to recall the definition of $f_{([t\alpha], \lambda)}$.
In this case,  $f_t$ rotates each of the $\lambda$--partition intervals $I_j$ by $t\alpha_j\lambda_j$ (mod $\lambda_j$). The nontrivial orbits 
of the action are the intervals $I_j$ for which the rotation rate $\alpha_j$ is nonzero. 

\begin{figure}[htb]   
\begin{center}
\includegraphics{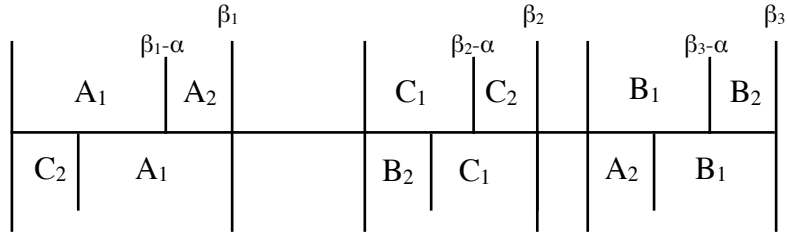}
\caption{A map conjugate to a restricted rotation \label{rotation conjugate} }
\end{center}
\end{figure}

 In general, a rotation subgroup $F = \{f_t\}$  
is conjugate in $\mathcal{E}$ to a subgroup $\{f_{([t\alpha], \lambda)} \}$ of a standard torus action; thus, $f_t = hf_{([t\alpha], \lambda)} h^{-1}$
for some $h \in \mathcal{E}$. The conjugacy $h$ may break
 a nontrivial $\{f_{([t\alpha], \lambda)}\}$\,--\,orbit into 
several intervals, but for each point $x$ in the interior of such an interval, the maps $f_t$ still locally translate $x$ by $t\alpha_j\lambda_j$. See an example 
in Figure~\ref{rotation conjugate}.
 
To make a precise statement, let $\text{Fix}(F) $
denote the set of global fixed points for a given one-parameter subgroup $F  = \{f_t\}$. Let $\mathcal{P}$ denote the 
set algebra of all finite unions of half-open intervals $[a, b)$ in $[0,1)$.

\begin{lemma} \label{local_rotation} A one-parameter subgroup $F = \{f_t\}$ of $\mathcal{E}$ is a rotation subgroup if 
and only if $\text{Fix}(F) \in \mathcal{P}$ and 
for all but finitely many $x\in[0,1)$, there exists $\alpha_x \in \mathbb{R}$ and $\epsilon_x > 0,$ such that
\[
f_t(x) = x + t\alpha_x,\ \  \text{if}\ |t| < \epsilon_x. 
\]  
\end{lemma}

\noindent  \emph {Proof}: It is easy to see that if $F$ is conjugate to a rotation subgroup, then the action of
$f_t$ satisfies the local condition stated in the lemma. In particular, if $f_t =  hf_{([t\alpha], \lambda)} h^{-1}$ for all $t$,
then the finite set of points that do not satisfy the condition is contained in the image under $h$ of the union of the set 
of discontinuities of $h$  and the set of partition interval endpoints induced by the length vector $\lambda$.

 Conversely, suppose $\text{Fix}(F)\in \mathcal{P}$ and $f_t$ is locally a rotation at all but finitely many points. 
Let $0=x_0 < x_1 < x_2 < \cdots < x_n = 1$ be the exceptional points, including all boundary points of $\text{Fix}(F)$.
 Over all $x\in (x_{i-1}, x_i)$ the
rotation speed $\alpha_x$ must be constant, since by definition it is locally constant. It remains to consider the
behavior of $f_t$ at the exceptional points.  

Consider the interval $I_j = \left[ x_{j-1}, x_j \right)$ of length $\lambda_j.$ Let $\alpha_j$ denote the constant rotation speed on the interior points of $I_j$. 
If $\alpha_j =0,$ then $I_j \subseteq \text{Fix}(F),$ so assume $\alpha_j \neq 0$.
By replacing $t$ with $-t$, it can be assumed that $\alpha_j$ is positive.
 For sufficiently small nonnegative $t$ and for sufficiently small $\delta$ such that $0 < \delta \ll \lambda_j$, the 
interval $(x_{j-1}, x_{j-1} + \delta)$ is translated by $t\alpha_j$ 
under $f_t.$  Since the maps $f_t$ are all right-continuous at $x_{j-1}$, it follows 
that $f_t(x_{j-1}) = x_{j-1} + t\alpha_j$ for all sufficiently small nonnegative $t$.

In fact, the group $f_t$ acts (locally) on all
of $(x_{j-1}, x_j)$ by translation by $t\alpha_j$; thus 
\[
 f_t(x_{j-1}) = x_{j-1} + t\alpha_j, \ \ \text{for } 0 \leq t < \frac{\lambda_j}{\alpha_j}. 
\] 
Consider what happens for $t = \lambda_j / \alpha_j$. First, suppose the interval $I_j$ is $f_t$--invariant. 
If $y = f_{(\lambda_j/\alpha_j)}(x_{j-1})$ is in the interior of $I_j$, then $[y, x_j)$ is a periodic orbit
properly contained in the orbit of $x_{j-1}$, which is impossible. Thus,
$f_{(\lambda_j/\alpha_j)}(x_{j-1}) = x_{j-1}$, and the action of $f_t$ on $I_j$ is 
globally a rotation action. 

In general, if $I_j$ is not invariant, suppose that $y$ is in the interior of some $I_k$, with $k\neq j,$
 since $y\in I_j$ would imply invariance. For small $t<0$, $f_t(y)$ is in $I_k$, since the $f_t$ locally act
as a rotation on the interior of $I_k$. However, $y = f_{(\lambda_j/\alpha_j)}(x_{j-1})$, and it is also
the case that $f_t(y)$ is in $I_j$ for small $t<0$, which is a contradiction. Thus,
 $f_{(\lambda_j/\alpha_j)}(x_{j-1})$ must be some other exceptional point $x_{k-1}$. The transformations $f_t$ all
 preserve Lebesgue measure, and by right-continuity $f_t(x_{k-1}) = x_{k-1} + t\alpha_k$ for small
 $t\geq 0.$ Thus, if $f_{(\lambda_j/\alpha_j)}(x_{j-1}) = x_{k-1}$, then $\alpha_j = \alpha_k$.
 Consequently, the orbit of $x_{j-1}$ is a finite union of intervals $I_k$, each of which has the same rotation speed.
 After applying a suitable conjugacy, each invariant collection of these intervals may be reassembled
 into a single invariant subinterval on which the conjugate action is
 a standard rotation action. \ \ $\square$    \\

 It is possible to improve on the previous lemma's recharacterization of rotation actions. In particular, since all interval 
 exchanges preserve Lebesgue measure on $\mathbb{T}$,
 the condition of a point $x$ having an orbit locally given by a rotation action can be weakened to the condition 
 that $t \mapsto f_t(x)$ is continuous for $t$ in a neighborhood of zero. If $x$ satisfies this weaker condition, 
 it is said to have a \emph{locally continuous orbit} under $f_t$.

 \begin{proposition} \label{local_cts_orbits_iff_rotn} A one-parameter subgroup $F = \{f_t\}$ is a rotation subgroup if and only if 
 $\text{Fix}(F) \in \mathcal{P}$ and 
 for all but finitely many $x \in [0,1)$, the function $\mathbb{R} \rightarrow [0,1)$ defined by
 \[  
 t \mapsto f_t(x)
 \]
 is continuous in some open neighborhood around $t=0$. 
 \end{proposition}

 \noindent \emph{Proof}: By applying the previous lemma, it suffices to show that if
  $x$ has a locally continuous
 orbit,   then there exists $\alpha_x$, such that $f_t(x) = x + t\alpha_x$ for all $t$ in a neighborhood of zero. 
 If $f_t(x) = x$ for all 
 $t$ in a neighborhood, then $x$ is a global fixed point of the action, and $\alpha_x = 0$ will suffice. 
 
 Suppose that $x$ is not a global fixed point, and 
 assume that the orbit $t \mapsto f_t(x)$  is continuous for $t\in [-\epsilon, \epsilon]$. 
  By reducing $\epsilon$ if necessary, it can be assumed that the function $t \mapsto f_t(x)$ is one-to-one
 on $[-\epsilon, \epsilon]$.

  By reversing the parameter $t$, it can be assumed that
  $t \mapsto f_t(x)$ is increasing on  $[-\epsilon, \epsilon]$. Define $\alpha \neq 0$ to satisfy 
 \[ 
     f_\epsilon(x) = x + \epsilon\alpha.
 \]             
 Next, for any $n\in \mathbb{N}$, consider the increasing sequence of points 
 \[
 f_{-\epsilon}(x), f_{{\frac{-(n-1)\epsilon}{n}} }(x), \ldots, x, f_{\frac{\epsilon}{n}} (x), f_{\frac{2\epsilon}{n}} (x),\, \ldots, f_\epsilon(x).
 \]
 Since all $f_t$ preserve Lebesgue measure $\mu$, 
 \[
     \mu\left(\left[ f_{\frac{(j-1)\epsilon}{n}}(x),\, f_{\frac{j\epsilon}{n}} (x)  \right) \right) =
      \mu\left(\left[ f_{\frac{(k-1)\epsilon}{n}}(x),\, f_{\frac{k\epsilon}{n}} (x)  \right) \right),
\]
for all $j$ and $k$ satisfying $-n < j,k \leq n$. Consequently, 
\[ 
   f_{\frac{j\epsilon}{n}} (x) = x +\left( \frac{j\epsilon}{n} \right)\!\alpha, \text{ for all integers $j$ such that } |j| \leq n.  
\]
Thus, $f_t(x) = x + t\alpha$ for a dense set of $t \in [-\epsilon,\epsilon]$, and by continuity of the orbit 
this holds at all $t \in [-\epsilon,\epsilon]$.   \ \ $\square$ \\

Therefore, to prove Theorem~\ref{one-parameter IET subgroups} it suffices to prove the following:

\begin{proposition} \label{R_action_has_loc_cts_orbits}
If F = $\{f_t\}$ is a continuous one-parameter subgroup of $\mathcal{E}$, then all but finitely many $x \in [0,1)$ have
locally continuous orbits and $\text{Fix}(F) \in \mathcal{P}$. \end{proposition}

Before giving a proof of this proposition, it will be shown that for any one-parameter subgroup $\{f_t\}$, the number of discontinuities
 of $f_t$ is bounded over all $t$ in some neighborhood of zero. Define the function $\delta\co \mathcal{E} \rightarrow \mathbb{N}$ by
\[
   \delta(f) =  n, \text{ where $f = f_{(\pi, \lambda)}$  for (unique) unpartitioned $\pi \in \Sigma_n, \lambda \in \Lambda_n$. } 
  \]
Equivalently, $\delta(f)$ returns the number of discontinuities of $f$, considered as a map $[0,1) \rightarrow [0,1)$, 
counting 0 as a discontinuity. 
 Note that $\delta(f) = n$ if and only if 
$f$ is in the image of the interior of the ($n$--1)--dimensional simplex $\overline{\Lambda}_{n}$ under the
parametrization $\Gamma_\pi$, for some unique unpartitioned $\pi \in \Sigma_{n}$. The simplex 
$\overline{\Lambda}_{n}$ is compact, and the number of unpartitioned permutations in $\Sigma_k$ for $k\leq n$
is finite. Since the parametrizations $\Gamma_\pi$ are continuous, the sets
\[
K_n = \left\lbrace f\in \mathcal{E}\co \delta(f) \leq n \right\rbrace
\]   
are compact.  Therefore, if $\delta(f) = n$, then $f$ is in the interior of $\mathcal{E} \setminus K_{n-1}$. 
Consequently, for all $g$ in some neighborhood of $f$, $\delta(g) \geq \delta(f)$;  in other words, the function $\delta$ is lower semicontinuous.

\begin{lemma} \label{delta(f_t) is bounded on compact subsets}
For any continuous one-parameter subgroup $F=f_t$, the function $t \mapsto \delta(f_t)$ is bounded 
on any compact subset of $\mathbb{R}$.
\end{lemma}

\noindent \emph{Proof}: Since $f_{s+t} = f_s\circ f_t$ for all $s,t \in \mathbb{R},$ it follows that
\begin{equation} \label{delta_ineq1}
\delta(f_{s+t}) \leq \delta(f_s) + \delta(f_t),\ \  s,t\in \mathbb{R}.
\end{equation}
This inequality records the fact that a composition of two interval exchange maps cannot have more discontinuities
than occur over both of its factors. From this inequality, it also follows that 
\begin{equation} \label{delta_ineq2}
\delta(f_{s+t}) \geq | \delta(f_s) - \delta(f_t)|,\ \ s,t \in \mathbb{R}.
\end{equation}

By \eqref{delta_ineq1}, if $\delta(f_t)$ is bounded for $t\in [-\epsilon, \epsilon]$, then 
$\delta(f_t)$ is bounded on all compact subsets. Thus, if $\delta(f_t)$ is unbounded on some compact 
subset, then $\delta(f_t)$ is
unbounded in any neighborhood of zero. In fact, the inequality \eqref{delta_ineq2} further implies that
$\delta(f_t)$ is unbounded in any neighborhood of any $t\in \mathbb{R}$. 

This local unboundedness and the semicontinuity of $\delta$ cannot coexist. To derive a contradiction, 
suppose that $\delta(f_t)$ is unbounded in any neighborhood of any $t$. Let
\[
  A_n = \left\lbrace t\in \mathbb{R}\co \delta(f_t) \leq n  \right\rbrace.
\] 
By the lower semicontinuity of $\delta$, the sets $A_n$ are closed, and their complements
\[
  B_n = \left\lbrace t\in\mathbb{R}\co \delta(t) > n \right\rbrace
\]  
are open. If $\delta$ is locally unbounded at every point, each set $B_n$ is dense in $\mathbb{R}$. 
However,
\[
    \bigcap B_n = \left\lbrace t\in \mathbb{R}\co \delta(f_t) > n
     \text{ for all } n \in\mathbb{N}  \right\rbrace = \varnothing,
\]
which is a contradiction by the Baire Category Theorem. Thus, $\delta(f_t)$ must be bounded 
on any compact subset of $\mathbb{R}$.  \ \ $\square $ \\

\noindent \emph{Proof of Proposition~\ref{R_action_has_loc_cts_orbits}}: Applying 
Lemma~\ref{delta(f_t) is bounded on compact subsets},
let \[n = \max\lbrace \delta(f_t)\co t \in [-1,1] \rbrace.\] By the lower semicontinuity of $\delta$, the set
$\lbrace t \in [-1,1]\co \delta(f_t)=n \rbrace$ is relatively open in $[-1,1]$. Therefore, there exists some $t_0 \in (-1,1)$
and $\epsilon >0$, such that $\delta(f_t) = n$ for $t \in (t_0 -\epsilon, t_0 + \epsilon)$. Let $\pi \in \Sigma_{n}$ be
the unique unpartitioned permutation such that $f_{t_0} \in \Gamma_\pi(\Lambda_{n})$. By Proposition~\ref{unique_coords}, the sets $\Gamma_\sigma(\Lambda_{n})$ are pairwise disjoint as $\sigma$ ranges over $\Sigma_{n}'$, the set of
unpartitioned permutations in $\Sigma_{n}.$ Since 
\[ f_{t_0} \in \mathcal{E} \setminus \left( \bigcup_{\sigma \in \Sigma_{n}' \setminus \{\pi\} } \Gamma_\sigma(\overline{\Lambda}_n)  \right) \]
and the sets $\Gamma_\sigma(\overline{\Lambda}_{n})$ are compact, it follows that $f_{t_0}$ is actually in the interior of the set denoted above.
 Consequently, after possibly replacing $\epsilon$ by a smaller value, it
follows that $f_t \in \Gamma_\pi(\Lambda_{n})$ for all $t \in (t_0 -\epsilon, t_0 + \epsilon)$.

In this situation, it can be seen that the paths 
\[
 t \mapsto f_t(x)
\]
are continuous in a neighborhood of $t_0$ for all but finitely many points, namely the discontinuity points of $f_{t_0}$.
For $t \in (t_0 -\epsilon, t_0 + \epsilon)$ let $\lambda^{(t)} \in \Lambda_{n}$ be such that
\[
 f_t = \Gamma_\pi\left(\lambda^{(t)}\right)\!,
 \] 
where the $\lambda^{(t)}$ vary continuously in $\Lambda_{n}$. Thus, if $x$ is an interior point of the interval $I_j$
induced by $f_{t_0}$, then 
\[ 
t \mapsto f_t(x) = x + \Omega_\pi\left(\lambda^{(t)}\right)_j 
\] 
is continuous in a neighborhood of $t_0$. Since $f_{-t_0}$ is continuous at all but a finite number of points,
the path $t \mapsto f_t(x)$ is continous in a neighborhood of zero for all but finitely many points.  

It remains to consider the set of global fixed points for $f_t$. As before, define $\beta_j^{(t)}$  in 
terms of $\lambda^{(t)}$ and let 
\[ I_j^{(t)} = \left[ \beta_{j-1}^{(t)}, \beta_j^{(t)}   \right).  \] 
Suppose the interior of $I_j^{(t_0)}$ contains a global fixed point $x$. Then for all $t$ in some
 $(t_0 -\epsilon, t_0 + \epsilon),$ the point $x$ is located in the interval $I_j^{(t)}$. Thus, for each $t$ in 
  $(t_0 -\epsilon, t_0 + \epsilon),$ the interval $I_j^{(t)}$ is fixed by $f_t$. In addition, the intervals $I_{j-1}^{(t)}$
  and $I_{j+1}^{(t)}$ cannot be fixed by $f_t$, since otherwise $\pi$ would be partitioned. As a result, the 
  boundary points $\beta_{j-1}^{(t)}$ and $\beta_j^{(t)}$ must be constant over $t \in  (t_0 -\epsilon, t_0 + \epsilon),$
  since otherwise there would be points fixed by $f_t$ for $t$ in some nonempty, proper open subset of $\mathbb{R}$, which
  is impossible. Thus, the set $\text{Fix}(F)$ of global fixed points for $f_t$ is a finite union of intervals 
  $I_j^{(t_0)}$, which implies that  $\text{Fix}(F)$ is a member of $\mathcal{P}$.  \ \ $\square$

%
%
%
\bibliographystyle{abbrv}

\bibliography{Novak_ref}

\end{document}